\documentclass[10pt,reqno]{amsart}

\usepackage{amssymb,amsmath,mathrsfs}

\usepackage[colorlinks=false,pdfborderstyle={/W 1},pdfpagelabels]{hyperref}

\def\bl{\begin{lemma}}
\def\el{\end{lemma}}
\def\bth{\begin{theorem}}
\def\eth{\end{theorem}}
\def\bc{\begin{corollary}}
\def\ec{\end{corollary}}
\def\bcj{\begin{conjecture}}
\def\ecj{\end{conjecture}}
\def\bpr{\begin{proposition}}
\def\epr{\end{proposition}}
\def\bde{\begin{definition}}
\def\ede{\end{definition}}
\newcommand{\be}{\begin{eqnarray}}
\newcommand{\ee}{\end{eqnarray}}

\def\E{\mathbb{E}}
\def\Pr{\mathbb{P}}

\def\QED{\hfill\qedsymbol}

\def\bdr{\partial}
\def\bdre{\partial_E}
\def\bdrv{\partial_V}
\def\bdriv{\partial_{iV}}
\def\bdrov{\partial_{oV}}
\def\ancip{\iota_\psi^*}   

\def\IPA#1{{\mathcal{IP}}^*_{#1}}

\def\Cl{\mathscr{C}}
\def\Blue{\mathscr{B}}
\def\Red{\mathcal{R}}
\def\PP{\mathscr{P}}
\def\QQ{\PP^*}
\def\UU{\mathscr{U}}
\def\W{{\mathcal W}}
\def\st{{\, : \ }}

\def\eps{\varepsilon}

\def\Z{{\mathbb Z}}

\def\Set{{\mathcal S}}
\def\A{{\mathcal A}}

\def\GG{{\mathcal G}}

\def\XX{{\mathcal X}}
\def\YY{{\mathcal Y}}

\newtheorem{theorem}{Theorem}[section]
\newtheorem{definition}{Definition}[section]
\newtheorem{lemma}[theorem]{Lemma}
\newtheorem{corollary}[theorem]{Corollary}
\newtheorem{proposition}[theorem]{Proposition}
\newtheorem{conjecture}[theorem]{Conjecture}

\theoremstyle{definition}
\numberwithin{equation}{section}

\input epsf.sty

\begin{document}

\title[Isoperimetric profile via exponential cluster repulsion]{A note on percolation on $\Z^d$:\\isoperimetric profile via exponential cluster repulsion}
\author{G\'abor Pete}
\address{{\it Address in 2008:} Microsoft Research, One Microsoft Way, Redmond, WA 98052-6399.\newline\indent
{\it Address in 2016:} Alfr\'ed R\'enyi Institute of Mathematics, Re\'altanoda u. 13-15, Budapest 1053 Hungary,\newline\indent
{\it and} Institute of Mathematics, Budapest University of Technology and Economics, Hungary}
\urladdr{http://www.math.bme.hu/~gabor}

\thanks{Partially supported by the Hungarian OTKA grant T049398.}

\date{Basically finished on April 24, 2008. Small correction in the proof of Theorem 1.2 on August 12, 2016.}

\begin{abstract}
We show that for all $p>p_c(\Z^d)$ percolation parameters, the probability that the cluster of the origin is finite but has at least $t$ vertices at distance one from the infinite cluster is exponentially small in $t$. We use this to give a short proof of the strongest version of the important fact that the isoperimetric profile of the infinite cluster basically coincides with the profile of the original lattice. This implies, e.g., that simple random walk on the largest cluster of a finite box $[-n,n]^d$ with high probability has $L^\infty$-mixing time $\Theta(n^2)$, and that the heat kernel (return probability) on the infinite cluster a.s.~decays like $p_n(o,o)=O(n^{-d/2})$. Versions of these results have been proven by Benjamini and Mossel (2003), Mathieu and Remy (2004), Barlow (2004) and Rau (2006). For general infinite graphs, we prove that anchored isoperimetric properties survive supercritical percolation, provided that the probability of the cluster of the origin being finite with large boundary decays rapidly; this is the case for a large class of graphs when $p$ is close to 1. As an application (with the help of some entropy inequalities), we give a short conceptual proof of a theorem of Angel, Benjamini, Berger and Peres (2006): the infinite percolation cluster of a wedge in $\Z^3$ is a.s.~transient whenever the wedge itself is transient.
\end{abstract}

\maketitle

\section{Introduction and results}\label{s.intro}

Isoperimetric inequalities on finite and infinite graphs are indispensable in studying the behavior of simple random walk (SRW) on the graph \cite{SC:StFlour,Woess,evolving}. Most importantly, a good isoperimetric profile implies fast mixing on a finite graph, or fast heat kernel decay on an infinite graph. It is important, from mathematical and physical points of view, to understand how robust these properties are under perturbations of the graph.
A standard question is as follows: consider supercritical Bernoulli$(p)$ edge-percolation on a transitive finite or infinite graph, then perform SRW on the giant or an infinite percolation cluster, respectively. {\it Do the most important properties of SRW survive percolation?} See the books \cite{Grimm,LPbook} for background on percolation.

On $\Z^d$ and its finite boxes, there is a large literature on this topic; the main results are the transience of the infinite cluster \cite{GKZ}, the right $d$-dimensional heat kernel decay and fast mixing \cite{BM:mixing,MaRe,Bar:rwperc}, and scaling to Brownian motion \cite{SiSz,BeBi:BM,MaPi}. For general transitive infinite graphs, the program was started by \cite{BLS:pertu}. For finite graphs other than boxes of $\Z^d$, only SRW on the giant component of the Erd\H os-R\'enyi random graph $G(n,p)$ has been understood fully \cite{BKW,evolution}. See \cite{anchsurv} for a recent survey on isoperimetry and SRW on percolation clusters.

In the Appendix of \cite{ChPePe}, our main discovery was that survival of the so-called anchored isoperimetry for infinite clusters can be deduced from an exponential decay of the probability that the cluster of the origin is finite but has a large boundary. This exponential decay has been proved only for large enough $p$ values; in fact, on $\Z^d$, when $d\geq 3$ and $p\in (p_c,1-p_c)$, only a stretched exponential decay holds. In the present note, we prove exponential decay on $\Z^d$, for all $p>p_c$, for a modified event, in which the boundary is not only large, but touches the infinite cluster at many places. Then, by refining a bit the main idea of \cite[Appendix]{ChPePe}, we prove survival of $d$-dimensional anchored isoperimetry. A good isoperimetic profile for the giant cluster of $[-n,n]^d$ will also follow, implying a strong mixing time result and $d$-dimensional heat kernel decay.
\medskip

In a connected bounded degree infinite graph $G(V,E)$, for $S\subset V$, the {\bf edge boundary} $\bdre S$ is the set of edges of $G$ with one endpoint in $S$, the other in $V\setminus S$. Similarly, the {\bf inner vertex boundary} $\bdriv S$ is the set of vertices that are in $S$ but have a neighbor outside $S$, while $\bdrov S:=\bdriv (V \setminus S)$ is the {\bf outer vertex boundary}. If it does not matter which boundary we are considering, we will drop the subscripts $E,V,i,o$. Furthermore, let $\overline{S}^G$ be the union of $S$ with all the finite connected components of $G \setminus S$; if $S$ is finite and connected, then so is this {\bf closure} $\overline{S}=\overline{S}^G$. The {\bf frontier} of $S$ is defined by $\bdr^+ S:=\bdr\overline{S}$, with the possible variations on $E,V,i,o$. For two percolation clusters, $\Cl_1$ and $\Cl_2$, a {\bf touching edge} is an edge of $G$ in $\bdre\Cl_1\cap\bdre\Cl_2$. The number of such edges will be denoted by $\tau(\Cl_1,\Cl_2)$. For supercritical percolation on $\Z^d$, the a.s.~unique infinite cluster is denoted by $\Cl_\infty$, while the cluster of the origin by $\Cl_o$. Our new percolation result is the following:

\bth\label{t.repulsion}
For $d\geq 2$ and any $p>p_c(\Z^d)$, there exists a $c_1=c_1(d,p)>0$ such that
\be\label{e.repul}
\Pr_p \Bigl( m\leq |\Cl_o| < \infty \hbox{ and }\, \tau(\Cl_o,\Cl_\infty)\geq t \Bigr)
\leq \exp\bigl(-c_1\max\{m^{1-1/d},t\}\bigr).
\ee
\eth

Setting $t=0$ in (\ref{e.repul}), the stretched exponential decay we get is a sharp classical result, due to Kesten and Zhang \cite{KZh} combined with the Grimmett-Marstrand theorem \cite{GrM}. Hence the exponential decay in $t$ is the novelty here. Nevertheless, our proof will be a modification of \cite{KZh}, so it naturally gives the $\exp(-c_1m^{1-1/d})$ part, as well. Moreover, (\ref{e.repul}) can probably be best understood from the perspective of \cite{KZh}. They prove the stretched exponential decay by showing that although for $p\in (p_c,1-p_c)$ the frontier $|\bdr^+\Cl_o|$ and the volume $|\Cl_o|$ are of the same order, $\Theta(m)$, there still exists a finite $N=N(p)$ such that the frontier of the set of vertices at distance at most $N$ from $\Cl_o$ is of size $\Theta(m^{1-1/d})$, and the probability of having such a large $\Cl_o$ is already exponential in this size. Therefore, having $\tau(\Cl_o,\Cl_\infty)\geq t \gg m^{1-1/d}$ means that $\Cl_\infty$ penetrates deep inside $\Cl_o$, going through tunnels of width less than $N$. As we will show, this has an exponentially small probability in $t$.

On nonamenable transitive graphs, there is conjecturally always an interval of $p$ values for which there are a.s.~infinitely many infinite clusters, see \cite{LPbook}. For this case, \cite{relentless} conjectured that no two infinite clusters can have infinitely many touching edges. This was recently proved by Tim\'ar \cite{Adam:nontouch} by an ingenious use of the Mass Transport Principle for unimodular transitive graphs (e.g.~all Cayley graphs). His argument might give some explicit decay for the probability that two neighboring vertices of an arbitrary unimodular transitive graph are in different clusters with at least $t$ touching edges, but getting the exponential decay rate in this general setting seems hard.

We use our Theorem \ref{t.repulsion} to prove the following sharp isoperimetric inequality:

\bth\label{t.ipabound} For $d\geq 2$ and $p>p_c(\Z^d)$, there are constants $\alpha(d,p)>0$ and $c_2(d,p)>0$ such that for the infinite cluster $\Cl_\infty=\Cl$, and for the edge frontier $\bdr_{\Cl}^+ S:=E_{\Cl}(\overline{S}^\Cl,\Cl\setminus\overline{S}^\Cl)$ inside $\Cl$,
\be\label{e.ipa}
\Pr_p\Bigl(\exists S\hbox{ connected}: o\in S\subset\Cl, M\leq |S| < \infty, \frac{|\bdr_{\Cl}^+ S|}{|S|^{1-1/d}}\leq \alpha\Bigr)\leq \exp\bigl(-c_2 M^{1-1/d}\bigr).
\ee
\eth

Considering only connected sets $S$ that contain a fixed origin $o$ is a natural restriction, since $\Cl_\infty$ has arbitrary large pieces with bad isoperimetry --- but they are typically far away from $o$. The following notion, introduced in \cite{Tho} and \cite{BLS:pertu}, is a general formulation of this idea. Take a connected bounded degree infinite graph $G(V,E)$, with a fixed $o \in V(G)$, and a positive function $\psi(\cdot)$ with $\lim_{x\to\infty}\psi(x)=\infty$. We say that $G$ satisfies an {\bf anchored $\psi$-isoperimetric inequality} if
\be
0< \ancip(G) := \lim_{n \to \infty}\, \inf \left\{{| \bdr S| \over \psi(|S|)}
\st o\in S\subset V(G),\,S \hbox{ is connected},\,n\leq
|S|<\infty \right\}.\label{e.ancpsi}
\ee
It is easy to see that the quantity $\ancip(G)$ does not depend on the choice of the basepoint $o$. The property $\ancip(G)>0$ is denoted by $\IPA{\psi}$, and, because of the bounded degrees, we can equally use $\bdr=\bdre$ or $\bdr=\bdrv$. For $\psi(x)=x$, this property is usually called {\bf anchored expansion} (or weak nonamenability), and for $\psi(x)=x^{1-1/d}$, we speak of {\bf $d$-dimensional anchored isoperimetry} $\IPA{d}$. Many probabilistic implications of isoperimetric inequalities remain true with this anchored version. Thomassen proved in \cite{Tho} that if $\IPA{\psi}$ holds with some function $\psi$ that satisfies
\be
\sum_{k=1}^\infty \psi(k)^{-2}<\infty,
\label{e.Tho}
\ee
then the graph contains a transient subtree, and so is transient itself. In particular, $\IPA{2+\eps}$ suffices for transience. Lyons, Morris and Schramm \cite{LMS} recently found a very nice few line proof of a refinement of Thomassen's result, resembling a converse to the Nash-Williams criterion; see also \cite{LPbook}. Vir\'ag proved in \cite{Vir:anch} the conjecture of \cite{BLS:pertu} that any bounded degree graph $G$ with anchored expansion has a non-amenable subgraph, and this subgraph is ``dense'' enough to ensure positive speed of SRW on $G$. On the other hand, it is not known if $\IPA{d}$ alone implies the usual $d$-dimensional heat kernel decay $p_n(o,o)=O(n^{-d/2})$. For more details and references see \cite{anchsurv}.

From Theorem \ref{t.ipabound}, the Borel-Cantelli lemma immediately implies that $\Cl_\infty$ a.s.~satisfies $\IPA{d}$. Moreover, we will also easily deduce the following isoperimetric profile:

\bc\label{c.profi} For all $p>p_c(\Z^d)$ there exist $c_3(d,p)>0$, $\alpha(d,p)>0$ and (for almost all percolation configurations $\omega$) an integer $N(\omega)$ such that for all $n>N(\omega)$, all connected subsets $S\subseteq \Cl_\infty \cap [-n,n]^d$ with size $|S|\geq c_3\,(\log n)^\frac{d}{d-1}$ have $|\bdr_{\Cl_\infty} S|\geq \alpha |S|^{1-1/d}$.
\ec

Conditioned on $o\in \Cl_\infty$, the walk on $\Cl_\infty$ started at $o$ cannot leave $[-n,n]^d$ in $n$ steps, so plugging this isoperimetric profile into the infinite graph heat kernel version of the Lov\'asz-Kannan bound \cite{LoKa}, proved by Morris and Peres \cite{evolving}, immediately gives that SRW on $\Cl_\infty$ has return probabilities $p_n(o,o)=O(n^{-d/2})$, for all $n>N(\omega)$. We will also prove the following finite version, which, in conjunction with the $L^\infty$-version of the Lov\'asz-Kannan bound, again from \cite{evolving}, implies that SRW on the largest cluster of $[-n,n]^d$ has $L^\infty$-mixing time $\Theta(n^2)$. The example of an infinite versus a finite depth regular tree shows that Corollary \ref{c.box} does not formally follow from Theorem \ref{t.ipabound}; nevertheless, the proofs of Theorems \ref{t.repulsion} and \ref{t.ipabound} can be modified to fit the finite case.

\bc\label{c.box} Let $\Cl$ be the largest cluster in percolation with $p>p_c(\Z^d)$ on the finite  box $[-n,n]^d$. Then $\exists\ c_3'(d,p)>0$ and $\alpha'(d,p)>0$ such that, with probability tending to 1, for all connected subsets $S\subseteq \Cl$ with size $c_3'\,(\log n)^\frac{d}{d-1}\leq |S|\leq |\Cl|/2$, we have $|\bdr_\Cl S|\geq \alpha' |S|^{1-1/d}$.
\ec

Corollary \ref{c.box} was first announced by Benjamini and Mossel \cite{BM:mixing}, but there were some gaps in their renormalization argument moving from $p$ values close to 1 to all $p>p_c(\Z^d)$. (These gaps seem repairable to us). A suboptimal bound $p_n(o,o)=O(n^{-d/2}(\log n)^{6d+14})$ was derived in \cite{HeHo}, while the true on-diagonal heat kernel and $L^\infty$-mixing time results were proved by Mathieu and Remy \cite{MaRe}. However, their isoperimetry results are weaker than ours. Barlow \cite{Bar:rwperc} proved the great result that a.s., for all large times $n>N_{x,y}(\omega)$, the heat kernel on $\Cl_\infty$ satisfies
$$
a_1n^{-d/2}\exp(-b_1\|x-y\|^2_1/n)\leq p_n(x,y)\leq a_2 n^{-d/2}\exp(-b_2\|x-y\|^2_1/n),
$$
with constants $a_i,b_i$ depending on $d$ and $p$, and random variables $N_{x,y}$ having at most a stretched exponential tail. Barlow did not state the sharp isoperimetric profile explicitly, but it can be deduced from his results (as shown to us by N.~Berger). Refining the approach of \cite{MaRe}, the preprint \cite{Rau} proves our Corollary \ref{c.profi} for $S\subseteq \Cl_\infty \cap [-n,n]^d$ with size $|S|\geq cn^\gamma$, arbitrary $c,\gamma>0$ and large enough $n$. Given the lengths of \cite{MaRe,Bar:rwperc,Rau}, we find our short proof of Theorem \ref{t.ipabound} and its corollaries very attractive. Independently, M.~Biskup has recently also constructed a short proof of Theorem \ref{t.ipabound} and Corollary \ref{c.profi}, along lines more similar to \cite{BM:mixing} than to our work. His proof appears in \cite{BBHK}, which paper shows that if the edges of $\Z^d$ are given i.i.d.~random conductances with a large tail at 0, then, in $d\geq 5$, the heat kernel decay is $\Theta(n^{-2})$; that is, the original decay $\Theta(n^{-d/2})$ does not survive this type of random perturbation. Such ``anomalous'' decay also happens when we move from supercritical to critical percolation: \cite{RWorientIIC} shows that the heat kernel decay on the incipient infinite cluster of oriented percolation on high dimensional $\Z^d$ is $\Theta(n^{-2/3})$. An analogous result for SRW on the critical Erd\H os-R\'enyi graph is proved in \cite{AsYu:diamix}.

Our proof of Theorem~\ref{t.repulsion} uses percolation renormalization, a method presently not available on other infinite graphs. However, for many graphs, for $p$ close to 1, it is easy to show a result even stronger than~(\ref{e.repul}), namely,
\be\label{e.peiexpon}
\Pr_p\bigl(|\Cl_o|<\infty,\ |\bdre^+ \Cl_o|=n\bigr)\leq \varrho^n
\ee
with $\varrho=\varrho(p)<1$ and all large $n$. This is the case, e.g., for Cayley graphs of finitely presented groups; see Theorem~\ref{t.kappa} below. The method of \cite[Appendix]{ChPePe} then implies the following:

\bth\label{t.general}
Suppose that $G$ satisfies $\IPA{\psi}$ with some $\psi\nearrow\infty$, and the exponential decay~(\ref{e.peiexpon}) holds for some $p$. Then $p$-a.s.~on the event that the open cluster $\Cl_o$ is infinite, $\Cl_o$ satisfies $\IPA{\psi}$.
\eth

As an application, we give a conceptual proof for a strengthening of the Grimmett-Kesten-Zhang theorem \cite{GKZ} of the transience of $\Cl_\infty$ in $\Z^3$:
survival of transience in more subtle situations. For an increasing positive function $h(\cdot)$, the {\bf wedge} $\W_h$ is the subgraph of $\Z^3$ induced by the vertices $V(\W_h)=\{(x,y,z) \st x \geq 0\ {\rm and}\ |z| \leq h(x)\}$. Terry Lyons \cite{TLy} proved that $\W_h$ is transient if{f}
\be
\sum_{j=1}^\infty {1\over jh(j)}<\infty\,.
\label{e.tlyons}
\ee
For example, (\ref{e.tlyons}) holds for $h(j)=\log^r j$ if{f} $r>1$. Now, the following holds.

\bth[Angel, Benjamini, Berger, Peres \cite{ABBP}]\label{t.wedge}
The unique infinite percolation cluster of a wedge $\W_h\subset\Z^3$ for any $p > p_c(\W_h)=p_c(\Z^3)$ is a.s.~transient if and only if $\W_h$ is transient.
\eth

The evolution of this result is that \cite{BPP:path} gave a new proof of \cite{GKZ}, and then, by sharpening those methods, H\"aggstr\"om and Mossel \cite{HM:low} verified the claim under the stronger condition $\sum_{j=1}^\infty {1\over j\sqrt{h(j)}}<\infty$, and asked whether Theorem~\ref{t.wedge} holds. We prove Theorem~\ref{t.wedge} under an additional mild concavity-type condition on $h(\cdot)$ that keeps technicalities to the minimum:

\bpr\label{p.LyonsThom}
If $h(\cdot)$ satisfies Lyons' condition (\ref{e.tlyons}) and there exists $\gamma>0$ such that $h(\delta x)\geq \gamma \delta h(x)$ for all $\delta\in [0,1]$, then $\W_h$ satisfies some $\IPA{\psi}$ with Thomassen's condition (\ref{e.Tho}).
\epr

The key step in the proof of this result is a projection type isoperimetric inequality in the wedge $\W_h$, similar to the Loomis-Whitney inequality \cite{Loomis}, which we show using some simple entropy inequalities. It was Han \cite{Han} and Shearer \cite{Shetal} who first proved entropy inequalities analogous to such isoperimetric inequalities, but it is unclear who noticed first that these are really the same results. See \cite{BalBol} for a concise treatment.

Given Proposition~\ref{p.LyonsThom}, our general Theorem~\ref{t.general} implies survival of transience for $p$ close to 1. Transience, unlike isoperimetry, is monotone w.r.t.~adding edges and vertices, so we can extend this result for all $p>p_c$ using a standard renormalization argument, and do not need a sophisticated result like Theorem~\ref{t.repulsion} showing the survival of isoperimetry itself for all $p>p_c$.
\medskip

\noindent{\bf Organization of paper.} Section~\ref{s.repuls} proves the percolation result Theorem~\ref{t.repulsion}. Section~\ref{s.isop} reaps its consequences to isoperimetry, Theorem~\ref{t.ipabound} and Corollaries~\ref{c.profi} and~\ref{c.box}. Section~\ref{s.general} shows that~(\ref{e.peiexpon}) holds for many graphs, and proves Theorem~\ref{t.general}. Finally, Section~\ref{s.wedge} deals with transient wedges.
\medskip

\noindent{\bf Some open problems.} There are many intriguing questions in the field. Does the giant cluster on the hypercube $\{0,1\}^n$ have mixing time polynomial in $n$? What is the heat kernel decay on the incipient infinite cluster of critical percolation on $\Z^d$? On an infinite transitive graph $G$, do transience, positive or zero speed, or certain heat kernel decay survive percolation for all $p>p_c(G)$? Does the analogue of our Theorem \ref{t.repulsion} hold on any transitive graph $G$? Does $\IPA{d}$ itself imply the heat kernel bound $p_n(o,o)\leq O(n^{-d/2})$? For a discussion of these and further questions, see \cite{anchsurv}.
\medskip

\noindent{\bf Acknowledgments.} I am grateful to Noam Berger, Marek Biskup and Yuval Peres for discussions and encouragement, to Russ Lyons for asking if my method in \cite[Appendix]{ChPePe} could work for transient wedges, and to Pierre Mathieu and \'Ad\'am Tim\'ar for comments on the manuscript. Also thanks to Yuval Peres and Perla Sousi for pointing out a small error in July 2016 in the proof of Theorem~\ref{t.ipabound}.

\section{Proof of the exponential cluster repulsion}\label{s.repuls}

We fix a positive integer $N$, whose $p$-dependent value will be determined later. We regard $N\Z^d$ as a graph naturally isomorphic to the lattice $\Z^d$, i.e., with adjacency relation $\|x-y\|_1 =N$. We will also use $N\Z^d_*$, the graph where adjacency is defined by $\|x-y\|_\infty =N$. We will use boxes of the form $B_{3N/4}(Nx):=\{y\in\Z^d: \|y-Nx\|_\infty\leq 3N/4 \}$, for $x\in\Z^d$. These will be called {\bf blocks}. The set of blocks will be thought of as the vertices of a graph naturally identified with $N\Z^d$.

The reason for considering two different adjacency relations are the following facts. While (a) is trivial from the definitions, the also quite innocent-looking (b) and (c) require careful proofs, which were executed in \cite[Lemma 2.1]{DPisz}. Recently, Tim\'ar \cite{Adam:conn} found a much simpler and more general proof. Recall the definitions of the closure $\overline S$ and the different boundaries $\bdr S$ from the Introduction.

\begin{itemize}{\it
\item[\bf(a)] For any finite $\Z^d$-connected set $A$, the vertex frontiers $\bdriv^+ A$ and $\bdrov^+ A$ are finite cutsets: any $\Z^d$-path connecting a vertex of $\overline A$ to a vertex of $\Z^d\setminus \overline A$ intersects both frontiers $\bdrv^+ A$.
\item[\bf(b)] The vertex frontiers $\bdrv^+ A$ need not be $\Z^d$-connected, but for $d\geq 2$ they are both $\Z^d_*$-connected.
\item[\bf(c)] Consider the finite box $B_n:=[-n,n]^d$, and a connected $A\subseteq B_n$. Let $A_i$ be the connected components of $B_n\setminus A$. Then all the vertex boundaries $\bdr_{V(B_n)} A_i$ are $\Z^d_*$-connected.}
\end{itemize}

As usual in renormalization, given a percolation configuration inside a block, we call the block {\bf good} if it has a cluster connecting all its $(d-1)$-dimensional faces, while all other clusters have diameter less than $N/5$. The basic result of static renormalization is that the probability that a given block is good tends to 1 as $N\to\infty$ \cite[Section 7.4]{Grimm}. Blocks not good will be called {\bf bad}.

From now on, we assume that $o\not\in\Cl_\infty$ and that the diameter of $\Cl_o$ is at least $N$. For any given cluster $\Cl$, a block $B$ is called {\bf $\Cl$-substantial} if $\Cl\cap B$ has a connected component of diameter at least $N/5$. The set of $\Cl$-substantial blocks will be denoted by $\Cl^N$; note that this is a connected subset of $N\Z^d$. Now, we color a block $B$ {\bf red} if it is $\Cl_o$-substantial but it has a neighbor that is not $\Cl_o$-substantial. In other words, the set of red blocks, denoted by $\Red$, equals $\bdriv (\Cl_o^N)$. Furthermore, we color a block {\bf blue} if it is both $\Cl_o$- and $\Cl_\infty$-substantial. The set of blue blocks is $\Blue$. Clearly, each pair of touching vertices is contained in at least one blue block, and in at most $2^d$. A block can be both red and blue. See Figure \ref{fig:RB}. Observe that a colored block is never good: on one hand, being blue implies the existence of two disjoint components of large diameter; on the other hand, in a good block $B$ that is $\Cl_o$-substantial, $\Cl_o$ must connect all the $(d-1)$-dimensional faces, which makes all the neighboring blocks $\Cl_o$-substantial, hence $B$ cannot be red. Our main claim is the following:

\begin{figure}[htbp]
\centerline{\epsfysize=2.8 true in \epsffile{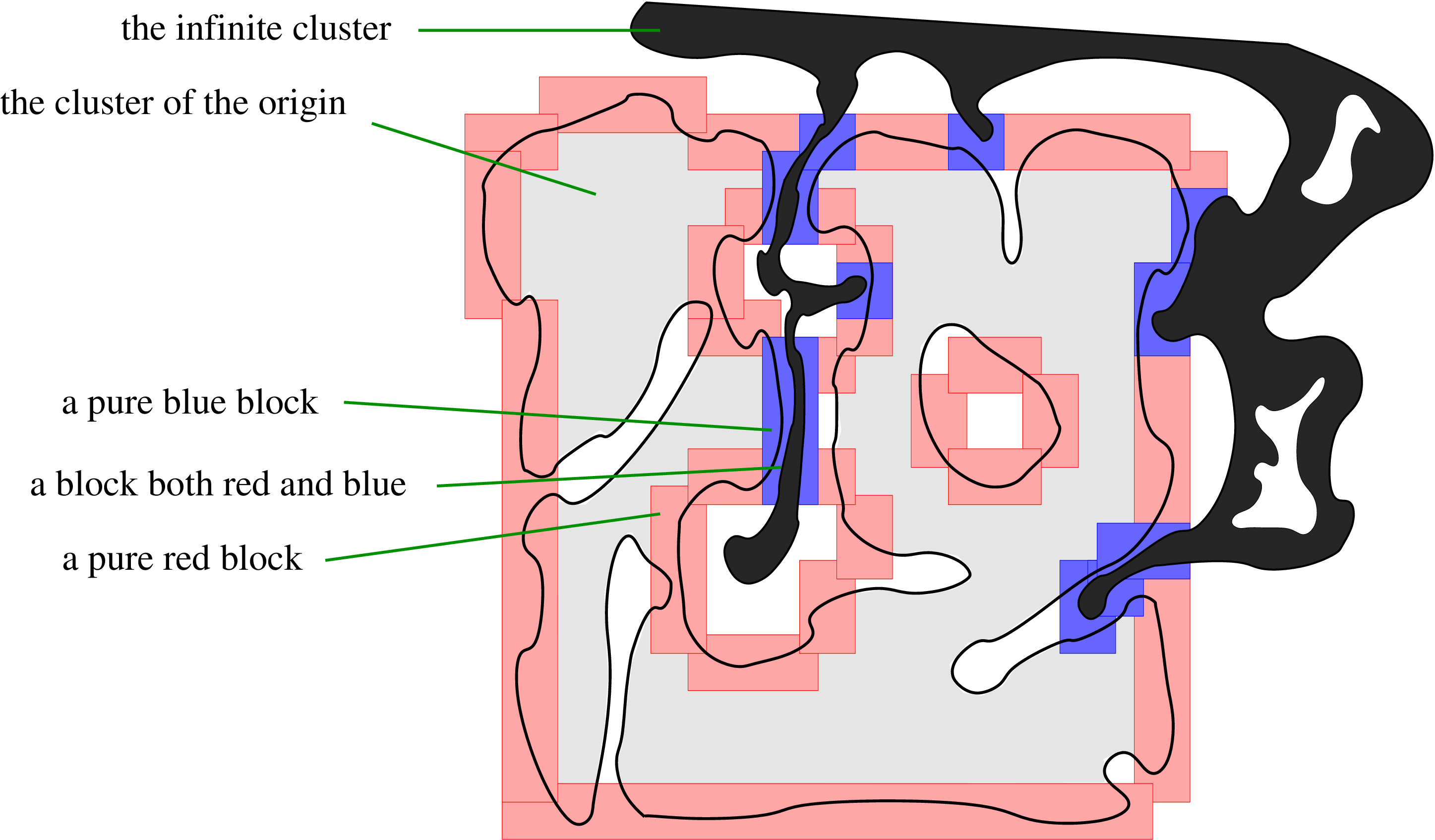}}
\caption{The clusters $\Cl_o$, $\Cl_\infty$ and the sets $\Red$, $\Blue$.
}\label{fig:RB}
\end{figure}


\bl\label{l.RBPQ}
On the event $\A_{m,t}:=\{|\Cl_o|=m \text{ and } \tau(\Cl_o,\Cl_\infty)\geq t\}$, the set $\Red \cup \Blue$ of colored blocks has an $N\Z^d_*$-connected subset of size $\geq c_4(N,d)\max\{m^{1-1/d},t\}$, contained in the box $B_m(o)$.
\el

\noindent{\it Proof.} We will first define $\PP$, a large $N\Z^d_*$-connected set that will be mostly colored. Then we will remove its uncolored parts and repair the resulting holes by adding some colored blocks, so that the augmented set $\QQ$ will be fully colored, $N\Z^d_*$-connected, and large.

Firstly, by Fact (b) above, the frontier $\bdriv^+(\Cl_o^N) \subseteq \Red$ is connected in $N\Z^d_*$, and is of size at least $c_5(d) m^{1-1/d}/N^d$, by the standard isoperimetric inequality in $\Z^d$. Secondly, consider $\overline{\Cl_o^N} \cap \Cl_\infty^N$. This set contains $\Blue$, whose size is at least $t/(2N)^d$. Now take the union
$$
\PP:=\bdriv^+(\Cl_o^N) \cup \bigl( \overline{\Cl_o^N} \cap \Cl_\infty^N\bigr).
$$
The set $\overline{\Cl_o^N} \cap \Cl_\infty^N$ can have several $N\Z^d$-connected components, but, by Fact (a), each component intersects $\bdriv^+(\Cl_o^N)$, which is an $N\Z^d_*$-connected set by (b). Therefore, the union $\PP$ is $N\Z^d_*$-connected, and its size is at least $\max\{c_5(d)m^{1-1/d}/N^d,t/(2N)^d\}$.

\begin{figure}[htbp]
\centerline{\epsfysize=1.8 true in \epsffile{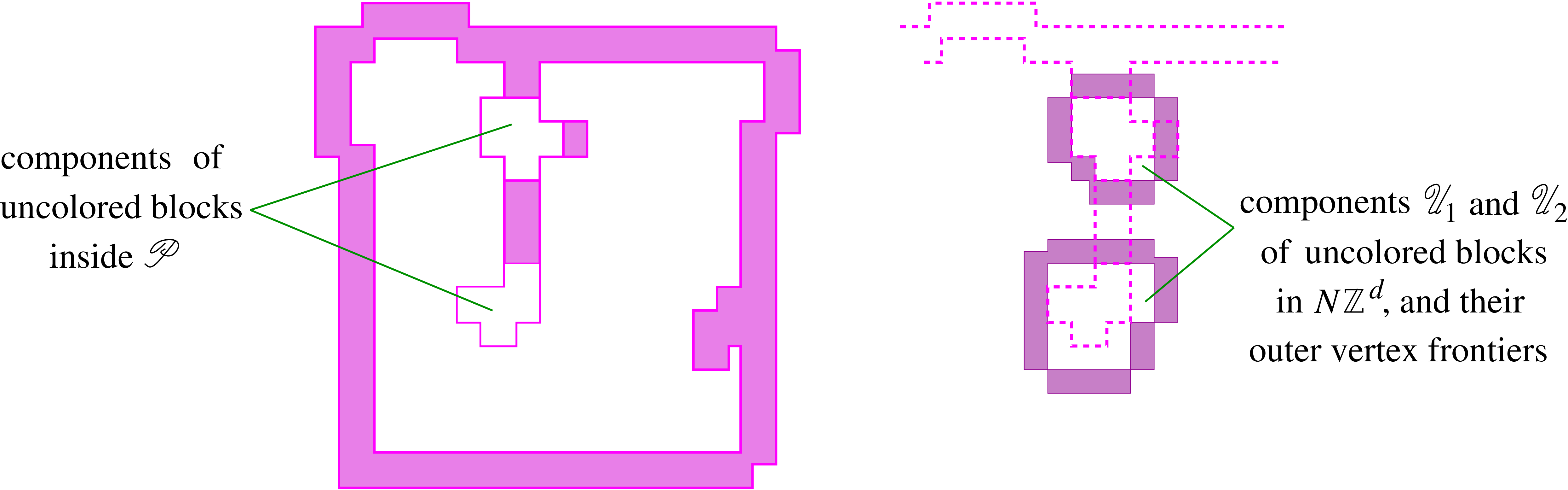}}
\caption{The $N\Z^d_*$-connected set $\PP$ copied from Figure \ref{fig:RB}.
}\label{fig:PP}
\end{figure}

However, $\PP$ may contain some uncolored blocks, listed as $U_1,\dots,U_k$. See the left side of Figure \ref{fig:PP}. The set of {\it all} uncolored blocks in $N\Z^d$ form connected components in $N\Z^d$: the infinite component $N\Z^d \setminus \overline{\Cl_o^N}$, and some finite ones, each separated from infinity by the red cutset $\bdriv^+(\Cl_o^N)$. Those finite components that contain at least one of the $U_i$'s will be listed as $\UU_1,\dots,\UU_\ell$. Clearly, each $U_i$ is in one of the $\UU_j$'s. We claim that
$$\QQ:=\Bigl(\PP\setminus \{U_j\}_{j=1}^k\Bigr) \cup \Bigl(\mathop{\bigcup}\limits_{j=1}^\ell \bdrov^+\UU_j\Bigr) =
\bigcup_{j=1}^{\ell} \Bigl((\PP\setminus \UU_j) \cup \bdrov^+\UU_j \Bigr)
$$ is an $N\Z^d_*$-connected subset of $\Red\cup\Blue$. See the right side of Figure \ref{fig:PP}.

For each $j$, the frontier $\bdrov^+\UU_j$ is a part of the boundary of an uncolored component, hence colored. Thus $\QQ \subseteq \Red\cup\Blue$; however, it is less clear that it is $N\Z^d_*$-connected. On the other hand, $\bdrov^+\UU_j$ is $N\Z^d_*$-connected by Fact (b), therefore $\PP^{**}:=\bigcup_{j=1}^{\ell} \Bigl((\PP\setminus \overline{\UU_j}) \cup \bdrov^+\UU_j \Bigr)$ is $N\Z^d_*$-connected. However, maybe it is smaller than $\QQ$. We will show in the next paragraph that $\PP\cap\overline{\UU_j}$ cannot have any colored blocks, which implies $\PP\cap\overline{\UU_j}=\PP\cap \UU_j$. Thus $\PP\setminus\overline{\UU_j}=\PP\setminus \UU_j$ and $\QQ=\PP^{**}$, which proves the claim that $\QQ$ is $N\Z^d_*$-connected.

Assume, on the contrary, that $\PP\cap\overline{\UU_j}$ has some colored block $D$. Since $\UU_j$ is disjoint from $\bdriv^+(\Cl_o^N)$, we have $\PP\cap\overline{\UU_j}\subseteq \overline{\Cl_o^N} \cap \Cl_\infty^N$, hence this colored $D$ is $\Cl_\infty$-substantial, and hence it is blue (it might at the same time be red). Thus there is a $\Cl_o^N$-path connecting $D$ to $\bdriv^+(\Cl_o^N)$, and a $\Cl_\infty^N$-path connecting $D$ to infinity. Both of these paths must intersect the uncolored cutset $\bdriv^+\UU_j$; let $F_\alpha\in\Cl_\alpha^N\cap \bdriv^+\UU_j$, where $\alpha\in\{o,\infty\}$. Since $F_\infty$ is uncolored but $\Cl_\infty$-substantial, it is non-$\Cl_o$-substantial. On the other hand, $F_o$ is $\Cl_o$-substantial. Both $F_\alpha$'s are in the $N\Z^d$-connected set $\UU_j$, hence, on any path connecting $F_o$ and $F_\infty$ inside $\UU_j$ there exists a $\Cl_o$-substantial block $F^*$ that has a non-$\Cl_o$-substantial neighbor. But then $F^*$ should be red by definition, while it is in $\UU_j$, hence uncolored --- a contradiction.

Thus we have $\QQ$, an $N\Z^d_*$-connected subset of $\Red\cup\Blue$, which contains all the colored blocks of $\PP$, hence its size is also at least $\max\{c_5(d)m^{1-1/d}/N^d,t/(2N)^d\}$. That it is inside the box $B_m(o)$ is clear from $|\Cl_o|=m$. Therefore, the lemma is proved.\QED
\medskip

\noindent{\it Proof of Theorem \ref{t.repulsion}}. Our Lemma \ref{l.RBPQ} says that $\A_{m,t}$ implies that the set of bad blocks contains an $N\Z^d_*$-connected subset of size at least $c_4(N,d)\max\{m^{1-1/d},t\}$ which is contained in $B_m(o)$. Whether a block is bad is independent of all the blocks which are not adjacent to it in $N\Z^d_*$. Therefore, a usual Peierls-argument for the graph $N\Z^d_*$ (see, e.g., \cite{KZh}, or Part (i) of Theorem \ref{t.kappa} below) gives that if the probability for a block to be bad is less than some $\eps>0$ whose value depends only on the graph structure $\Z^d_*$, then the probability of $\A_{m,t}$ is less than $\exp\bigl(-c_1(N,d)\max\{m^{1-1/d},t\}\bigr)$. As mentioned above, the point of renormalization is exactly that the probability of a block being bad is less than any $\eps>0$ if $N$ is large enough, thus the proof of Theorem \ref{t.repulsion} is complete.\QED

\section{Proof of the $d$-dimensional isoperimetric inequalities}\label{s.isop}

\noindent{\it Proof of Theorem \ref{t.ipabound}}. We will denote the infinite percolation cluster $\Cl_\infty$ simply by $\Cl$. For a connected subgraph $S\subseteq \Cl \subseteq \Z^d$, denote $\tilde\bdr_\Cl S:=E_{\Z^d}(\overline{S}^\Cl,\Cl \setminus \overline{S}^\Cl)$, the set of edges in $\Z^d$ with one endpoint in $S$, the other in the unique infinite component of $\Cl\setminus S$. Note that $\tilde\bdr_\Cl S \cap E(\Cl) = \bdr_\Cl^+ S$, i.e., $\bdr_\Cl^+ S$ is the set of edges in $\tilde\bdr_\Cl S$ that are open and hence belong to $\Cl$. Consider now the events
$$
\XX(m,t,s):=\{\exists S\text{ connected}: o\in S\subset \Cl,\ |S|=m, |\tilde\bdr_\Cl S|=t, |\bdr_\Cl^+ S|=s\}.
$$
Our goal is to bound from above the quantity
\be\label{e.rosszak}
\Pr_p\Bigl(\exists S\text{ conn.}: o\in S\subset\Cl, |S|\geq M, \frac{|\bdr_\Cl^+ S|}{|S|^{1-1/d}}\leq \alpha\Bigr) = \sum_{m\geq M}\sum_{s=1}^{\lfloor \alpha m^{1-1/d}\rfloor}\sum_{t=s}^{2md} \Pr_p(\XX(m,t,s)),
\ee
where we used that the number of edges leaving $S$ is at most $2d|S|$. For the events
$$
\YY(m,t):=\{|\Cl_o|=m \quad\text{and}\quad \tau(\Cl_o,\Cl_\infty)=t\},
$$
our Theorem \ref{t.repulsion} says that, for some $c=c(d,p)>0$,
\be\label{e.Ybound}
\Pr_p(\YY(m,t))\leq \exp\bigl(-c\max\{m^{1-1/d},t\}\bigr).
\ee
Given a configuration $\omega\in\XX(m,t,s)$ and a corresponding set $S\ni o$, define a new configuration $F(\omega,S)\in\YY(m,t)$ by redeclaring the edges in $\bdr^+_\Cl S$ to be closed. For a given $\omega'\in\YY(m,t)$, there are ${t\choose s}$ pre-images $(\omega,S)$ under $F$. For each $\omega\in\XX$ there is at least one $S$, hence, writing $Q=\frac{p}{1-p}>1$,
\be\label{e.XY}
\Pr_p(\XX(m,t,s)) \leq {t\choose s} Q^s \, \Pr_p(\YY(m,t)).
\ee
Combining (\ref{e.Ybound}) and (\ref{e.XY}) will give an upper bound on (\ref{e.rosszak}). The summations over $s$ and $t$ in (\ref{e.rosszak}) can be rewritten, for any $K>1$ that we will fix soon, as
$$
\sum_{t=1}^{\lfloor K \alpha m^{1-1/d}\rfloor} \; \sum_{s=1}^{\min\{t,\lfloor \alpha m^{1-1/d} \rfloor\}}\ \Pr_p(\XX(m,t,s)) \ + \sum_{t=\lfloor K \alpha m^{1-1/d}\rfloor+1}^m \sum_{s=1}^{\lfloor \alpha m^{1-1/d}\rfloor} \Pr_p(\XX(m,t,s)) \ =:\ \Set_1+\Set_2.
$$
We have $\Set_1\leq \sum_{t=1}^{\lfloor K \alpha m^{1-1/d}\rfloor} (1+Q)^t\exp\bigl(-c\max\{m^{1-1/d},t\}\bigr)$. For $\alpha=\alpha(d,p,K)$ small enough, this is at most $\exp\bigl(-(c/2) m^{1-1/d}\bigr)$. To bound $\Set_2$, we are using the straightforward estimate
\be
\sum_{s=1}^{\alpha n} {\beta n \choose s} Q^s \leq \exp\Bigl(\alpha\bigl(1+\log (\beta/\alpha)+\log Q\bigr)\,n\Bigr)\qquad \text{for }\beta\geq\alpha,
\ee
applied with $n=m^{1-1/d}$ and $t=\beta n$, where $\beta > K \alpha$. If $K=K(d,p)$ is large enough compared to how large $Q$ 
 and how small $c$ are, then, for any $\alpha,\beta>0$ with $\beta > K \alpha$,
$$\alpha(1+\log(\beta/\alpha)+\log Q) < \alpha(1+2\log(\beta/\alpha))  < (c/2) \beta\,.$$
Hence
$\Set_2 \leq \sum_{t=\lfloor K \alpha m^{1-1/d}\rfloor+1}^m \exp\bigl(-(c/2)t\bigr)$. Putting together our bounds on $\Set_1$ and $\Set_2$, we get that (\ref {e.rosszak}) is at most $\exp\bigl(-c' M^{1-1/d}\bigr)$ for some $c'=c'(d,p)$.\QED
\medskip

\noindent{\it Remark.} The decay rate in (\ref{e.ipa}) is sharp for a simple reason. Take $r\in\Z^+$, an edge $e_r \in \bdre [-r,r]^d$, some $0< \rho < \Pr_p(o\in \Cl_\infty)$, and consider the event $\A_r:=\Bigl\{o\in\Cl_\infty$, $e_r\in E(\Cl_\infty)$, $|\Cl_\infty \cap [-r,r]^d| > \rho r^d\Bigr\}$.
Then $\Pr_p(\A_r)>c>0$. For $\omega\in\A_r$, define $\hat\omega$ by redeclaring all of $\bdre [-r,r]^d \cap E(\Cl_\infty)$ but $e_r$ to be closed. By counting preimages and the cost of redeclaration, the set $\hat\A_r$ of $\hat\omega$'s has probability at least $\exp(-Cr^{d-1})$, and on $\hat\A_r$, the connected set $S:=\Cl_\infty \cap [-r,r]^d$ has $|S| >\rho r^d$ but $|\bdr^+_\Cl S|=1$.
\medskip

\noindent{\it Proof of Corollary \ref{c.profi}}. Consider percolation on the infinite lattice $\Z^d$. Then, by Theorem \ref{t.ipabound} and a union bound, $\Pr_p\Bigl(\exists x\in [-n,n]^d \hbox{ and }S\hbox{ connected}: x\in S\subset \Cl_\infty,\  |S|\geq M, \frac{|\bdr_{\Cl}^+ S|}{|S|^{1-1/d}}\leq \alpha\Bigr)$ is at most $(2n)^d\exp\bigl(-c_2 M^{1-1/d}\bigr)$. If $M\geq c_3(\log n)^{\frac{d}{d-1}}$ with $c_3>(d+2)/c_2(d,p)$, then this probability is at most $O(1/n^2)$, so the Borel-Cantelli lemma finishes the proof.\QED
\medskip

\noindent{\it Proof of Corollary \ref{c.box}}. We first need a finite box version of Theorem \ref{t.repulsion}. For this, consider two disjoint clusters $\Cl_1,\Cl_2\subseteq B_n=[-n,n]^d$ with $m\leq |\Cl_1|\leq |\Cl_2|$ and $\tau(\Cl_1,\Cl_2)\geq t$.

For simplicity, we assume that $n$ is divisible by $N$, and define the red set $\Red$ using $\Cl_1$, and the blue set $\Blue$ using $\Cl_1$ and $\Cl_2$, analogously to what we did in Section \ref{s.repuls}. We have $\Red,\Blue\subseteq NB_{n/N}$. Now let $\Delta:=\bdriv^{+\Cl_2}(\Cl_1)$ be the set of vertices in $\Cl_1$ that have a neighbor in the connected component of $B_n \setminus \Cl_1$ that contains $\Cl_2$. Consider $\Delta^N \subseteq NB_{n/N}$. By Fact (c) above, $\Delta^N$ is an $N\Z^d_*$-connected subset of $\Red$, and, similarly to Fact (a), it is easy to see that each $N\Z^d$-connected component of $\overline{\Cl_1^N} \cap \Cl_2^N$ intersects $\Delta^N$.
Furthermore, by the finite box isoperimetric inequalities of \cite{BolL} or \cite{DPisz}, the size of $\Delta^N$ is at least $c(N,d) m^{1-1/d}$. (This is the step that works for $\Z^d$ but would break down on a finite ball of a regular tree.) Therefore, our new
$$\PP:=\Delta^N \cup \bigl(\overline{\Cl_1^N} \cap \Cl_2^N\bigr)$$
is again an $N\Z^d_*$-connected set of size at least $c(N,d)\max\{m^{1-1/d},t\}$. Exactly as before, the corresponding set $\QQ$ is a large $N\Z^d_*$-connected subset of $\Red \cup \Blue$. Finally, renormalization gives that for $p>p_c(\Z^d)$, the probability of having disjoint clusters $\Cl_1,\Cl_2$ with $m\leq |\Cl_1|\leq |\Cl_2|$ and $\tau(\Cl_1,\Cl_2)\geq t$ is at most $C n^d \exp\bigl(-c_1(N,d)\max\{m^{1-1/d},t\}\bigr)$.

Plugging this into the proof of Theorem \ref{t.ipabound}, we get: with probability going to 1, for any subset $S$ of the giant cluster $\Cl$ such that both $S$ and $\Cl\setminus S$ are connected, and $c_3'(d,p)\,(\log n)^\frac{d}{d-1}\leq |S|\leq |\Cl\setminus S|$, we have $\frac{|\bdr_\Cl S|}{|S|^{1-1/d}}\geq \alpha'(d,p)$. We have to extend this result to all connected subsets $S$ that are large enough. This is exactly the content of the not very hard Lemma 2.6 of \cite{BM:mixing}, and we are done.\QED
\medskip

\section{Survival of anchored isoperimetry on general graphs}\label{s.general}

Consider a bounded degree infinite graph $G$, with a fixed vertex $o$. Let $q_n(G)$ be the number of minimal edge cutsets of cardinality $n$ separating $o$ from infinity. Assume that
\be
\kappa(G):=\limsup_{n\to\infty} q_n(G)^{1/n}<\infty,\label{e.kap}
\ee
which quantity does not depend on the basepoint $o$. (\ref{e.kap}) is known to be satisfied in many situations: when $G$ is the Cayley graph of a finitely presented group that is not a finite extension of $\Z$, or is quasi-isometric to such a Cayley graph \cite{BaB,Adam:cut}; when $G$ is a planar graph with polynomial growth and isoperimetric dimension bigger than 1 \cite{PPP}; when $G$ has anchored expansion \cite{ChPePe}. We will see that transient wedges of $\Z^3$ also satisfy  (\ref{e.kap}). And why is this useful for us?

\bth\label{t.kappa} Consider edge-percolation on a bounded degree infinite graph $G(V,E)$.
\item{\bf (i)} If $\kappa(G)<\infty$, then $p_c(G)\leq 1-1/\kappa(G)$, and the exponential decay~(\ref{e.peiexpon}) holds for $p>1-1/\kappa(G)$.
\item{\bf (ii)} Suppose that $G$ satisfies $\IPA{\psi}$ with some $\psi\nearrow\infty$, and that the exponential decay (\ref{e.peiexpon}) holds for some $p$. Then $p$-almost surely on the event that the open cluster $\Cl_o$ is infinite, $\Cl_o$ satisfies $\IPA{\psi}$.
\eth

\noindent{\it Proof}. Part (i) is a standard Peierls argument. For any $p>1-1/\kappa(G)$ fixed, let $\eps \in (0,\,1/\kappa-1+p)$, and $N$ is so large that $q_n(G)<(\kappa+\eps)^n$ for all $n>N$. Then the expected number of minimal edge cutsets of size $n$ with all edges being closed is $q_n(G)(1-p)^n < (1-\eps^2)^n$. This expectation is an upper bound on the probability of having any such cutset, hence (\ref{e.peiexpon}) is proved. Moreover, if $N$ is large enough, then the probability of having any closed cutset of size larger than $N$ is strictly less than 1. Now, $\kappa<\infty$ easily implies the existence of some integer $r$ such that the ball of radius $r$ around $o$ has $|\bdre^+ B_r(o)|>N$. With positive probability, this ball is not separated from infinity. But the event $\{B_r(o)\subseteq \Cl_o\}$ is independent from this separation, and it has positive probability, so both events together occur with positive probability, and then $\Cl_o$ is infinite.

Part (ii) can be proved following the Appendix of \cite{ChPePe} almost verbatim. In the language of our above proof of Theorem \ref{t.ipabound}, the argument is as follows. Consider the events
\be
\XX(m,s)&:=&\{|\Cl_o|=\infty, \hbox{ and }\exists S\hbox{ connected}: o\in S\subset \Cl,\ |\bdre^+S|=m, |\bdr_\Cl^+ S|=s\},\cr
\YY(m)&:=&\{|\Cl_o|<\infty,\ |\bdre^+ \Cl_o|=m\}.\nonumber
\ee
Now, for $\omega\in\XX(m,s)$ with $s\leq\alpha m$, and a corresponding $S\ni o$, we redeclare the edges in $\bdr^+_\Cl S$ to be closed, and get $F(\omega,S)=\omega'\in\YY(m)$. Therefore,
\be\label{e.XYmegint}
\sum_{s\leq\alpha m}\Pr_p(\XX(m,s)) \leq \sum_{s\leq\alpha m} {m\choose s} Q^s \Pr_p(\YY(m)).
\ee
To bound $\Pr_p(\YY(m))$ from above, we use (\ref{e.peiexpon}) in place of (\ref{e.Ybound}). On the other hand, for any $\eps>0$, if $\alpha$ is small enough, then $\alpha m {m\choose \alpha m} Q^{\alpha m} < (1+\eps)^m$. Thus we get an exponential decay for (\ref{e.XYmegint}), and the Borel-Cantelli lemma gives a positive lower bound on the ratios $|\bdr^+_\Cl S|/|\bdr^+_E S|$. Now, for an arbitrary connected set $o\in S\subset \Cl$, we can take its closure $\overline{S}^\Cl$ inside the graph $\Cl$. It is easy to see that $|\bdr^+_\Cl \bigl(\overline{S}^\Cl\bigr)|/|\bdre^+ \bigl(\overline{S}^\Cl\bigr)| \leq |\bdr_\Cl S|/|\bdre S|$, which implies that  $\IPA{\psi}$ survives. \QED
\medskip

\section{Percolation on transient wedges}\label{s.wedge}

\noindent{\it Proof of Proposition \ref{p.LyonsThom}}. Consider a connected subset $S$ in $\W_h$, containing the origin $o$, of volume $|S|=v$ and boundary $|\bdre S|=w$. We are going to show that there exist a constant $c=c(\W_h)\in (0,1)$ and some $k=k(S) \in \Z^+$ such that
\be
w \geq c\sqrt{h(k) v}
\label{e.egy}
\ee
and
\be
w \geq v/k\,.
\label{e.ketto}
\ee
We claim that these imply
$$w \geq c\sqrt{vf(v)}\,,\qquad \hbox{where}\qquad f(v):=h\Big(\sqrt{v/h(\sqrt{v})}\Big)\,.$$
Indeed, if $k \geq \sqrt{v/h(\sqrt{v})}$, then the monotonicity of $h$ implies $h(k)v\geq f(v)v$, and the claim follows from (\ref{e.egy}). If $k \leq \sqrt{v/h(\sqrt{v})}$, then $v/k \geq \sqrt{vh(\sqrt{v})}\geq \sqrt{vf(v)}$, so the claim follows from (\ref{e.ketto}).

That is, $\W_h$ satisfies $\IPA{\psi}$ with $\psi(v)=\sqrt{vf(v)}$. So the last thing we need for~(\ref{e.Tho}) is that $\sum_{k=1}^\infty (vf(v))^{-1}<\infty$.
Our condition $h(\delta x)\geq \gamma\delta h(x)$ implies $h(\sqrt{v})\leq C \sqrt{v}$ with $C=h(1)/\gamma$. These and the monotonicity of $h$ give
$({\gamma}/{\sqrt{C}})\, h(v^{1/4})\leq h\big({v^{1/4}}/{\sqrt{C}}\big)\leq h\big(\sqrt{v/h(\sqrt{v})}\big)$.
On the other hand, by a change of variables, (\ref{e.tlyons}) implies $\sum_{v=1}^{\infty} (vh(v^{1/4}))^{-1}<\infty$, and we are done.

\medskip
\noindent{\it Remark.} Note here that the seemingly natural choice $\psi(v):=\sqrt{v h(v)}$ does not work, as can be easily checked, e.g., for $h(x)=x^\alpha$, any $\alpha\in (0,1]$.
\medskip

We still need to prove (\ref{e.egy}) and (\ref{e.ketto}). For this, we will use some simple entropy inequalities. For random variables $\xi,\eta$ with values in a finite set $A$, their {\bf entropy} and {\bf conditional entropy} are
$$H(\xi):=\sum_{a\in A}\Pr(\xi=a)\log \frac{1}{\Pr(\xi=a)}\qquad\hbox{and}\qquad H(\xi\mid\eta):=H(\xi,\eta)-H(\eta)\,,$$
where $H(\xi,\eta)$ is the entropy of the r.v.~$(\xi,\eta)$. We will be using two basic inequalities: entropy is maximized by the uniform measure, i.e., $H(\xi)\leq \log |A|$, and $H(\xi\mid\eta,\,\zeta) \leq H(\xi\mid\eta)$.

Consider the projections $P_x,P_y,P_z$ in the three coordinate directions of $\Z^3$. Given $S \subset \W_h$, we use the notation $S(x,y,\cdot):=S\cap (x,y,\Z)$ for the sections of $S$. Let $w_x:=|P_x(S)|$, $w_y:=|P_y(S)|$, and $w_z:=\big|\big\{(x,y)\in P_z(S): |S(x,y,\cdot)| < |\W_h(x,y,\cdot)|=2h(x)+1\big\}\big|$. Note that $w \geq  w_x + 2w_y +w_z$.

Let $(X,Y,Z)$ be a uniform random point of $S$. Note that $H(X,Y,Z) = \log v$, while $H(Y,Z)\leq \log w_x$ and $H(X,Z)\leq \log w_y$. On the other hand, from the basic properties of conditional entropy, one easily gets $H(Y,Z)+H(X,Z) \geq H(X,Y,Z) + H(Z)$. This gives
\be\label{e.wH}
w_x\, w_y \geq v \exp(H(Z))\,.
\ee
Now we let $k := v/(w_x+w_z)$. In the proof below, this quantity will play the role of some sort of weighted average size of $S$ in the $x$ direction. Since $w\geq w_x+w_z$, we obviously have~(\ref{e.ketto}). The key step now will be to show that
\be\label{e.Hhk}
H(Z) \geq \log \big(c'\, h(k)\big)
\ee
for some $c'=c'(\W_h)>0$, because then~(\ref{e.wH}) and the inequality between the geometric and arithmetic means imply $w_x+w_y \geq 4\sqrt{v\, c'\, h(k)}$, and then~(\ref{e.egy}) follows immediately.

Decompose $S$ into $S_\text{full}:=\{(x,y,z)\in S: |S(x,y,\cdot)|=2h(x)+1\}|$ and $S_\text{miss}:=S\setminus S_\text{full}$, with sizes $v_\text{full}+v_\text{miss}=v$.
Denote $k_\text{full}:=|S_\text{full}|/|P_x(S_\text{full})|\geq v_\text{full}/w_x$.
If $(Y_\text{full},Z_\text{full})$ is picked uniformly in $P_x(S_\text{full})$, and $\xi_\text{full}:= |S_\text{full}(\cdot,Y_\text{full},Z_\text{full})|$, then $\E\xi_\text{full}=k_\text{full}$, while the r.v.~$|S_\text{full}(\cdot,Y,Z)|$ conditioned on $(X,Y,Z) \in S_\text{full}$ is the size-biased version of $\xi_\text{full}$. Therefore,
$$
\Pr\Big( |S_\text{full}(\cdot,Y,Z)| \leq \eps k_\text{full} \;\Big|\; (X,Y,Z)\in S_\text{full}\Big)
= \sum_{j\leq \eps k_\text{full}} \frac{j\Pr(\xi_\text{full}=j)}{\E\xi_\text{full}}
\leq \eps \Pr(\xi_\text{full}\leq \eps k_\text{full}) \leq \eps\,,
$$
for any $\eps>0$. It follows immediately that
\be\label{e.full}
\Pr\big( X \leq \eps k_\text{full} \;\big|\; (X,Y,Z)\in S_\text{full}\big) \leq \eps\,.
\ee

Similarly, if we let $h_\text{miss}:=|S_\text{miss}|/|P_z(S_\text{miss})| = v_\text{miss}/w_z$, then
\be\label{e.miss}
\Pr\Big( |S_\text{miss}(X,Y,\cdot)| \leq \eps h_\text{miss} \;\Big|\; (X,Y,Z)\in S_\text{miss}\Big)
\leq \eps\,.
\ee

Denoting $\nu:=v_\text{full}/v$ and $\rho:=w_x/(w_x+w_z)$, we have $k_\text{full} = (\nu/\rho) \, k$ and $h_\text{miss}=(1-\nu)/(1-\rho)\,k$. Introducing the r.v.~$\zeta:=|S(X,Y,\cdot)|$, the inequalities (\ref{e.full}) and (\ref{e.miss}) translate to
$$
\Pr\big(\zeta \leq h(\delta k) \mid (X,Y,Z)\in S_\text{full}\big) \leq \frac{\rho}{\nu}\,\delta\,,
\qquad\hbox{and}\qquad
\Pr\big(\zeta \leq \delta k \mid (X,Y,Z)\in S_\text{miss}\big) \leq \frac{1-\rho}{1-\nu}\,\delta\,,
$$
for any $\delta>0$. Since $h(\delta k)\leq \delta k$, these together give
\be\label{e.zhdk}
\Pr\big(\zeta \leq h(\delta k)\big)\leq \delta\,.
\ee
Finally, notice that this concentration result and our condition $h(\delta x)\geq \gamma \delta h(x)$ imply
\be
H(Z)\geq H(Z \mid X,Y)=\E(\log h(\zeta)) & \geq & \sum_{j\geq 1} \frac{1}{2^j} \log h(k/2^j)
\nonumber\\
& \geq &\sum_{j\geq 1}\frac{\log \big(\gamma h(k) \big)}{2^j} - \sum_{j\geq 1}\frac{j \log 2}{2^j} \geq \log h(k) - C\,, \nonumber
\ee
and (\ref{e.Hhk}) is proved.
\QED
\medskip

\noindent{\it Proof of Theorem \ref{t.wedge} with the extra assumption on $h(\cdot)$}. Because of the translation invariance in the $y$ direction and the amenability of $\W_h$, we can have only one infinite cluster a.s., at any $p$; see \cite{LPbook}. The fact that $p_c(\W_h)=p_c(\Z^3)$ whenever~(\ref{e.tlyons}) holds will be clear from what follows.

One direction is standard: if $\W_h$ is recurrent, then any subgraph of it is also such, by Rayleigh's monotonicity principle \cite{LPbook}. Conversely, when $\W_h$ is transient: from Proposition~\ref{p.LyonsThom} we now that $\W_h$ has $\IPA{\psi}$ with some $\psi$ satisfying Thomassen's condition~(\ref{e.Tho}). Furthermore, (\ref{e.kap}) holds, by the following argument. Firstly, let $\GG$ be the graph whose vertices are the edges of $\W_h$, and two such edges are adjacent in $\GG$ if they have some endpoints that are $\Z^3_*$-adjacent. Each degree in $\GG$ is at most a constant $D$, and any minimal edge-cutset in $\W_h$ separating $o$ from infinity is a connected subgraph of $\GG$. Secondly, if the distance in $\W_h$ of an edge-cutset from the origin is at least $t$, then its cardinality is at least $t$, because of its intersection with the plane $(x,y,0)\subset \W_h$. Altogether, the number of edge-cutsets of size $n$ is at most $O(n^3)\Delta^n$, hence $\kappa(\W_h)<\infty$, indeed. Theorem~\ref{t.kappa} now gives that the infinite cluster at $p>1-1/\kappa(\W_h)$ also satisfies $\IPA{\psi}$, and thus it is transient.

Now we want to extend this result for any $p>p_c(\Z^3)$; this will be almost the same as in~\cite{ABBP}. Recall the definitions of a {\bf block}, a {\bf good} block and a {\bf $\Cl$-substantial} block from Section~ \ref{s.repuls}, w.r.t.~an integer $N=N(p)$. Let $\W_h(N)$ be the set of blocks that are contained in $\W_h$; we will think of $\W_h(N)$ as a subgraph of $N\Z^3\simeq \Z^3$ or $N\Z^3_*\simeq \Z^3_*$. The monotonicity of $h$ implies that $\W_h(N)$ is infinite and connected for any $N$. Again, $\W_h(N)$ has at most one infinite cluster at any $p$, and $\W_h(N)$ satisfies the same $\IPA{\psi}$ as $\W_h$.

In Section \ref{s.repuls} we used that the probability for a block to be good is at least $1-\eps$ for $N$ large. A stronger statement is the Antal-Pisztora renormalization lemma \cite[Proposition 2.1]{chemical}, which also follows from the general Liggett-Schonmann-Stacey domination theorem \cite{domination}. Applied to $\W_h$, it says that for all $p>p_c(\Z^3)$ and $\eps>0$ there is an $N$ so that the process $\tilde\Pr_{p,N}$ of good blocks stochastically dominates Bernoulli$(1-\eps)$ percolation $\Pr_{1-\eps,N}$ on $\W_h(N)$, and the $\Cl_\infty$-substantial blocks form a unique infinite component on $\W_h(N)_*$, denoted by $\Cl_\infty(N)$. Moreover, the a.s.~transience of $\Cl_\infty$ on $\W_h$ under $\Pr_p$ would follow from the a.s.~transience of $\Cl_\infty(N)$ under $\tilde\Pr_{p,N}$.

It is not difficult to see that the cutset exponent~(\ref{e.kap}) satisfies
$\kappa((\W_h(N)_*)=\kappa((\W_h)_*)<\infty$, for any fixed $N$ and $h$. The same holds for the vertex cutset exponent $\kappa_V$ that can be defined analogously. Given $p>p_c(\Z^3)$, take $N$ so large that $\tilde\Pr_{p,N}$ dominates $\Pr_{1-\eps,N}$ site percolation on $\W_h(N)$, where $1-1/\kappa_V((\W_h)_*) < 1-\eps < 1$. Then, by Theorem~\ref{t.kappa}, we have a unique infinite Bernoulli$(1-\eps)$-cluster on $\W_h(N)_*$, which is transient if $h$ satisfies~(\ref{e.tlyons}). Rayleigh's monotonicity principle implies that $\Cl_\infty(N)$ is also transient, so, finally, $\Cl_\infty$ is such, too.\QED


\medskip

\begin{thebibliography}{ABBP06}

\bibitem[ABBP06]{ABBP} O. Angel, I. Benjamini, N. Berger and Y. Peres.
Transience of percolation clusters on wedges. {\it Elec. J. Probab.} {\bf 11} (2006),
655--669. {\tt arXiv:math.PR/0206130}

\bibitem[AntP96]{chemical} P. Antal and A. Pisztora. On the chemical distance in supercritical Bernoulli percolation. {\it Ann. Probab.} {\bf 24} (1996), 1036--1048.

\bibitem[BabB99]{BaB} E. Babson and I. Benjamini. Cut sets and normed cohomology, with applications to percolation. {\it Proc. Amer. Math. Soc.} {\bf 127} (1999), 589--597.

\bibitem[BalBo]{BalBol} P. Balister and B. Bollob\'as.
Projections, entropy and sumsets.
{\it Preprint}, {\tt arXiv:0711.1151v1 [math.CO]}

\bibitem[Bar04]{Bar:rwperc} M. T. Barlow. Random walks on supercritical
percolation clusters. {\it Ann. of Probab.} {\bf 32} (2004), 3024--3084.

\bibitem[BarJKS]{RWorientIIC} M. T. Barlow, A. J\'arai, T. Kumagai and G. Slade.
Random walk on the incipient infinite cluster for oriented percolation in high dimensions.
{\it Comm. Math. Phys.}, to appear. {\tt arXiv:math.PR/0608164}

\bibitem[BeKW]{BKW} I. Benjamini, G. Kozma and N. Wormald. The mixing time of the giant component of a random graph. {\it Preprint}, {\tt arXiv:math.PR/0610459}.

\bibitem[BeLS99]{BLS:pertu} I. Benjamini, R. Lyons and O. Schramm.
Percolation perturbations in potential theory and random
walks, {\it In: Random walks and discrete potential theory (Cortona,
1997)}, Sympos. Math. XXXIX, M. Picardello and W. Woess (eds.),
Cambridge Univ. Press, Cambridge, 1999, pp. 56--84. {\tt arXiv:math.PR/9804010}

\bibitem[BeM03]{BM:mixing} I. Benjamini and E. Mossel. On the mixing time of simple random walk on the super critical percolation cluster. {\it Probab. Theory Related Fields} {\bf 125} (2003), no. 3, 408--420. {\tt arXiv:math.PR/0011092}

\bibitem[BePP98]{BPP:path} I. Benjamini, R. Pemantle and Y. Peres. Unpredictable paths and percolation. {\it Ann. Probab.} {\bf 26} (1998), 1198--1211.

\bibitem[BerB07]{BeBi:BM}
N. Berger and M. Biskup. Quenched invariance principle for simple random walk on percolation clusters. {\it Probab. Theory Related Fields} {\bf 137} (2007), 83--120. {\tt arXiv:math.PR/0503576}

\bibitem[BerBHK]{BBHK}
N. Berger, M. Biskup, C. Hoffman and G. Kozma.
Anomalous heat kernel decay for random walk among bounded random conductances.
{\it  Ann. Inst. H. Poincar\'e}, to appear. {\tt arXiv:math.PR/0611666}

\bibitem[BoL91]{BolL} B. Bollob\'as and I. Leader. Edge-isoperimetric inequalities in the grid.
{\it Combinatorica} {\bf 11} (1991), 299--314.

\bibitem[ChPP04]{ChPePe} D. Chen and Y. Peres, with an appendix by G. Pete.
Anchored expansion, percolation and speed. {\it Ann. Probab.} {\bf 32}
(2004), 2978--2995. {\tt arXiv:math.PR/0303321}

\bibitem[ChGFS86]{Shetal}
F.~R.~K. Chung, R.~L. Graham, P. Frankl and J.~B. Shearer.
Some intersection theorems for ordered sets and graphs.
{\it J. Combinatorial Theory A} {\bf 43} (1986), 23--37.

\bibitem[DeP96]{DPisz} J-D. Deuschel and A. Pisztora. Surface order large deviations for high-density percolation. {\it Probab. Th. Rel. Fields} {\bf 104} (1996), no.~4, 467--482.

\bibitem[FoRe]{evolution} N. Fountoulakis and B. Reed.
The evolution of the mixing rate. {\it Preprint}, {\tt arXiv:math.CO/0701474}.

\bibitem[Gri99]{Grimm}
G. Grimmett. {\it Percolation. Second edition}. Grundlehren der
Mathematischen Wissenschaften, 321. Springer-Verlag, Berlin, 1999.

\bibitem[GKZ93]{GKZ} G. Grimmett, H. Kesten and Y. Zhang. Random walk on
the infinite cluster of the percolation model. {\it Probab. Th. Rel.
Fields}, {\bf 96} (1993), 33--44.

\bibitem[GrM90]{GrM} G. Grimmett and J. Marstrand. The supercritical phase of percolation is well-behaved. {\it Proc. Roy. Soc. London Ser. A} {\bf 430} (1990), 439--457.

\bibitem[H\"aPS99]{relentless} O. H\"aggstr\"om, Y. Peres and R.~H. Schonmann.
Percolation on transitive graphs as a coalescent process: Relentless merging followed by simultaneous uniqueness. {\it In: Perplexing Problems in Probability (M. Bramson and R. Durrett, ed.)}, pages 69--90, Boston, Birkh\"auser. Festschrift in honor of Harry Kesten.

\bibitem[H\"aM98]{HM:low} O. H\"aggstr\"om and E. Mossel. Nearest-neighbor walks with low predictability profile and percolation in $2+\epsilon$ dimensions.
{\it Ann. Probab.} {\bf 26} (1998), 1212--1231.

\bibitem[Han78]{Han} T. S. Han.
Nonnegative entropy measures of multivariate symmetric correlations.
{\it Information and Control} {\bf 36} (1978), 133--156.

\bibitem[HeH05]{HeHo} D. Heicklen and C. Hoffman. Return probabilities of a simple random walk on
percolation clusters. {\it Electron. J. Probab.} {\bf 10} (2005), 250--302.

\bibitem[KeZh90]{KZh} H. Kesten and Y. Zhang. The probability of a large finite cluster in supercritical Bernoulli percolation. {\it Ann. Probab.} {\bf 18} (1990), 537--555.

\bibitem[Kozm]{PPP} G. Kozma. Percolation, perimetry, planarity. {\it Rev. Math. Iberoam.} {\bf 23} (2007), no.~2, 671--676. {\tt arXiv:math.PR/0509235}.

\bibitem[LiSS97]{domination} T.~M. Liggett, R.~H. Schonmann and A.~M. Stacey.
Domination by product measures. {\it Ann. Probab.} {\bf 25} (1997), 71--95.

\bibitem[LoWh49]{Loomis} L.H. Loomis and H. Whitney.
An inequality related to the isoperimetric inequality.
{\it Bull. Amer. Math. Soc.} {\bf 55} (1949), 961--962.

\bibitem[LoKa99]{LoKa} L. Lov\'asz and R. Kannan. Faster mixing via average conductance. {\it Proceedings of the 27th Annual ACM Symposium on theory of computing}, 1999.

\bibitem[LyMS]{LMS} R. Lyons, B. Morris and O. Schramm. Ends in uniform spanning forests.
{\it Preprint,} {\tt arXiv:0706.0358 [math.PR]}.

\bibitem[LyPer]{LPbook} R. Lyons, with Y. Peres. {\it Probability on
trees and networks}. Book in preparation, present version is at {\tt
http://mypage.iu.edu/\~{}rdlyons}.

\bibitem[LyT83]{TLy} T. Lyons. A simple criterion for transience of a reversible Markov chain. {\it Ann. Probab.} {\bf 11} (1983), 393--402.

\bibitem[MaP07]{MaPi} P. Mathieu and A.~L. Piatnitski. Quenched invariance principles for random walks on percolation clusters. {\it Proc. R. Soc. Lond. Ser. A Math. Phys. Eng. Sci.} {\bf 463} (2007), 2287--2307. {\tt arXiv:math.PR/0505672}

\bibitem[MaR04]{MaRe} P. Mathieu and E. Remy. Isoperimetry and heat kernel decay on percolation clusters. {\it Ann. Probab.} {\bf 32} (2004), 100--128. {\tt arXiv:math.PR/0301213}

\bibitem[MoP05]{evolving} B. Morris and Y. Peres. Evolving sets, mixing and heat kernel bounds. {\it Prob. Th. Rel. Fields} {\bf 133} (2005), no. 2, 245--266. {\tt arXiv:math.PR/0305349}

\bibitem[NaPer]{AsYu:diamix} A. Nachmias and Y. Peres. Critical random graphs: diameter and mixing time. {\it Ann. Prob.}, to appear. {\tt arXiv:math.PR/0701316}

\bibitem[Pete]{anchsurv} G. Pete. Anchored isoperimetry, random walks, and  percolation: a survey with open problems. {\it In preparation}.

\bibitem[Rau]{Rau} C. Rau. Sur le nombre de points visit\'es par une marche al\'eatoire sur un amas infini de percolation. {\it Preprint}, {\tt arXiv:math.PR/0605056}.

\bibitem[SaC97]{SC:StFlour} L. Saloff-Coste. Lectures on finite Markov chains. {\it In: Lectures on probability theory and statistics (Saint-Flour, 1996)}, pages 301--413. Lecture Notes in Math., 1665, Springer, Berlin, 1997.

\bibitem[SiSz04]{SiSz} V. Sidoravicius and A-S. Sznitman.
Quenched invariance principles for walks on clusters of percolation or among random conductances. {\it Probab. Theory Related Fields} {\bf 129} (2004), 219--244.

\bibitem[Tho92]{Tho} C. Thomassen. Isoperimetric inequalities
and transient random walks on graphs. {\it Ann. Probab.} {\bf 20} (1992),
1592--1600.

\bibitem[Tim06]{Adam:nontouch} \'A. Tim\'ar. Neighboring clusters in Bernoulli percolation. {\it Ann. Probab.} {\bf 34} (2006), no.~6. 2332--2343. {\tt arXiv:math.PR/0702873}

\bibitem[Tim07]{Adam:cut} \'A. Tim\'ar. Cutsets in infinite graphs. {\it Combin. Probab. \& Comput.} {\bf 16} (2007), no.~1, 159--166.

\bibitem[Tim]{Adam:conn}  \'A. Tim\'ar. Some short proofs for connectedness of boundaries. {\it Preprint}, {\tt http://www.math.ubc.ca/$\sim$timar}.

\bibitem[Vir00]{Vir:anch} B. Vir\'ag. Anchored expansion and
random walk. {\it Geom. Funct. Anal.} {\bf 10} (2000), no. 6, 1588--1605.
{\tt arXiv:math.PR/0102199}

\bibitem[Woe00]{Woess} W. Woess.
{\it Random walks on infinite graphs and groups}.
Cambridge Tracts in Mathematics Vol.~138, Cambridge University Press, 2000.

\end{thebibliography}
\end{document}